\documentclass[12pt]{article}

\usepackage{amsmath}
\usepackage{amssymb}
\usepackage{amsfonts}
\usepackage{amsthm}
\usepackage{graphicx}
\usepackage{setspace}
\usepackage{float}
\usepackage{booktabs}
\usepackage{threeparttable}
\usepackage{multirow}
\usepackage{enumerate}
\usepackage{natbib}
\setcitestyle{aysep={}}
\usepackage{url}
\usepackage[hidelinks]{hyperref}
\usepackage{xcolor}
\usepackage{bm}
\usepackage{indentfirst}

\newif\ifblinded
\blindedfalse

\theoremstyle{plain}

\newtheorem{theorem}{Theorem}

\newtheorem{example}{Example}
\newtheorem{assumption}{Assumption}
\newtheorem{corollary}{Corollary}

\DeclareMathOperator{\Cov}{Cov}
\DeclareMathOperator{\Var}{Var}
\DeclareMathOperator{\Corr}{Corr}

\allowdisplaybreaks

\begin{document}

\def\spacingset#1{\renewcommand{\baselinestretch}%
{#1}\small\normalsize} \spacingset{1}

\ifblinded
{
  \bigskip
  \bigskip
  \bigskip
  \begin{center}
    {\LARGE\bf Glivenko--Cantelli Theorems for Integrated Volatility Functionals in Pure-Jump Semimartingales with an Application to Cryptocurrency Markets}
  \end{center}
  \medskip
}
\else
{
  \title{\bf Glivenko--Cantelli Theorems for Integrated Volatility Functionals in Pure-Jump Semimartingales with an Application to Cryptocurrency Markets}
  \author{
    Dachuan Chen\thanks{School of Economics, Singapore Management University, Singapore.}
    \quad Jia Li\textsuperscript{*}
    \quad Tian Xie\textsuperscript{*}
    \quad and Chengxin Yan\textsuperscript{*}}
  \maketitle
}
\fi

\begin{abstract}
\noindent We develop a two-step procedure for estimating integrated volatility functionals, defined through the occupation measure of the latent spot volatility process, when the asset price is a pure-jump semimartingale. In the first step, block-based estimators formed from absolute powers of high-frequency increments uniformly approximate local averages of powers of volatility. In the second step, these estimates are aggregated into an empirical occupation measure. Since price increments have infinite variance in this setting, arguments based on local Gaussianity are unavailable, and the uniform theory instead rests on maximal inequalities tailored to the stable regime. We establish Glivenko--Cantelli-type uniform consistency over classes of bounded monotone, Lipschitz-in-parameter, and locally H\"older test functions. These results deliver consistent estimation of volatility occupation times and quantiles, together with an argmax-consistency theory for $M$-estimators built on nonparametrically recovered latent processes. We further propose a stability-based rule for selecting the power index 
of the volatility estimator, which tracks an infeasible ex ante optimal 
choice closely in Monte Carlo experiments.
 An application to high-frequency cryptocurrency markets illustrates the framework in a jump-dominated, heavy-tailed environment.
\end{abstract}

\vfill

\noindent%
 \vspace{0.5cm}

\newpage
\spacingset{1.75}

\section{Introduction}
\label{sec:intro}

Volatility is a central object in the statistical analysis of financial data \citep{engle2004risk,andersen2003modeling}. Many quantities of interest in the study of volatility can be expressed as an integrated volatility functional of the form
\begin{align}
\label{IVF}
	\int_0^T g(\sigma_s)\,ds,
\end{align}
where $\sigma_t$ denotes the latent spot volatility of the asset price process and $g$ is a test function. Suitable choices of $g$ deliver a broad range of risk-related quantities. Power functions yield the integrated variance and quarticity \citep{jacod2013quarticity}, indicator functions yield volatility occupation times and their quantiles \citep{li2013volatility}, exponential functions yield the volatility Laplace transform \citep{todorov2012realized}, and richer families of test functions generate integrated moment conditions for the estimation of parametric volatility models \citep{li2016generalized}. The functional \eqref{IVF} may equivalently be viewed as the action on $g$ of the occupation measure induced by the volatility path. Estimating \eqref{IVF} uniformly over a class of test functions is therefore a Glivenko--Cantelli-type problem for an empirical occupation measure, and this is the perspective adopted in the present paper.

For It\^{o} semimartingales with a nondegenerate Brownian component, inference for \eqref{IVF} based on discretely sampled prices is by now well developed. One strand of the literature establishes consistency and central limit theory for smooth test functions \citep{barndorff2004econometric,mykland2009inference,jacod2012discretization,jacod2013quarticity,li2016inference,renault2017efficient}, with subsequent refinements concerning efficiency \citep{jacod2019estimating,li2019efficient,li2021efficient} and plug-in estimation of continuous-time regression models \citep{li2017adaptive}. Another strand studies discontinuous test functions, notably the volatility occupation time and its quantiles \citep{li2013volatility,li2016estimating}, and \citet{li2021glivenko} establish Glivenko--Cantelli theorems for the empirical occupation measure over general classes of test functions. These results rely, in an essential way, on the local Gaussianity of the price increments induced by the Brownian driver. By contrast, when the price process is of pure-jump type, existing results are available only for specific parametric families of test functions, namely power functions \citep{woerner2007inference,todorov2011limit} and exponential functions \citep{todorov2012realizeda}. To the best of our knowledge, a theory of uniform estimation of \eqref{IVF} over general classes of test functions in the pure-jump setting is not yet available. Developing such a theory is the goal of this paper.\footnote{\citet{yan2026nonparametric} study nonparametric estimation and inference for the spot volatility process itself in the same pure-jump setting. That paper does not consider integrated volatility functionals or the associated Glivenko--Cantelli-type uniform estimation problem, which are the focus of the present paper.}

The pure-jump setting is relevant well beyond any single market. Stable processes and related heavy-tailed models have a long tradition in statistics and applied probability, with applications ranging from signal processing and network traffic to insurance and physics \citep{samorodnitsky1994stable,nolan2020univariate}. Pure-jump models also have a long history in financial econometrics, where L\'evy processes of normal inverse Gaussian type provide parsimonious models for asset returns \citep{barndorff1997processes}, and the L\'evy-driven Ornstein--Uhlenbeck framework of \citet{barndorff2001non} models volatility itself as a pure-jump process. In high-frequency finance, the prevalence of jumps is commonly measured by the jump-activity index $\beta\in(0,2]$, where the value $\beta=2$ corresponds to the Brownian benchmark and lower values indicate dynamics dominated by large jumps and heavier tails \citep{ait2009estimating,todorov2015jump,kolokolov2022estimating}. Using this framework, \citet{todorov2011volatility} document that high-frequency data on the VIX index behave as a pure-jump process of infinite variation. Cryptocurrency markets provide further evidence, as high-frequency crypto returns display pronounced discontinuities and heavy tails \citep{liu2021risks,liu2022common,scaillet2020high}, and \citet{kolokolov2022estimating} finds that Bitcoin, in contrast to conventional exchange rates, is well described by a pure-jump process without a Brownian component. Our own estimates for Bitcoin (BTC) and the TRUMP token (TRUMP) support this evidence. We accordingly develop the estimation theory under a pure-jump specification with $\beta\in(1,2)$, and we revisit the cryptocurrency data in Section~\ref{sec:empirical} as an empirical illustration.

We develop a semiparametric two-step procedure for estimating \eqref{IVF}. In the first step, we construct a block-based spot volatility estimator from the $p$th absolute powers of high-frequency increments, with $p\in(0,\beta)$, and we show that it uniformly approximates the corresponding local averages of volatility across all blocks, provided that the number of blocks does not grow too quickly (Theorem~\ref{th:spot-uniform}). The main technical difficulty is that the absence of a Brownian component replaces Gaussian scaling with $\beta$-stable scaling. Increments then have infinite variance, so maximal inequalities based on second moments are unavailable, and the uniform control of the estimation error across a growing number of blocks requires arguments tailored to the stable regime.

In the second step, we plug the first-step estimates into the occupation measure to form an empirical occupation measure. We first establish pointwise consistency for general, possibly random, test functions under mild continuity and boundedness conditions (Theorem~\ref{th:OM-point-consist}). We then prove uniform consistency over three classes of test functions. For the class of bounded monotone functions, which contains all indicator functions, we obtain a Glivenko--Cantelli theorem that delivers uniformly consistent estimation of volatility occupation times, together with consistent estimation of the associated volatility quantiles (Theorem~\ref{th:OM-uniform-consist-G_M} and its corollaries). For a Lipschitz-in-parameter class, we obtain uniform convergence of plug-in criterion functions over compact parameter spaces, which yields an argmax-consistency theory for $M$-estimators based on nonparametric spot volatility estimates (Theorems~\ref{th:OM-uniform-consist-G_L} and~\ref{th:M-estimators}). For a locally H\"older class, we cover the nonlinear transformations, such as power and exponential functions, that are commonly used in the analysis of volatility functionals (Theorem~\ref{th:OM-uniform-consist-G_H}). Taken together, these results extend the Glivenko--Cantelli-type theory of \citet{li2021glivenko} from the Brownian semimartingale setting to the pure-jump setting.

Implementing the procedure requires choosing the power index $p$, and the Brownian default $p=2$ is no longer appropriate in the pure-jump setting. We propose a simple selection rule based on a stability consideration. Since the target functional does not depend on $p$, the rule searches for a value of $p$ at which the estimate is insensitive to small perturbations of $p$. The rule is easy to implement, requires only a rough preliminary estimate of $\beta$, and delivers robust and adequate finite-sample performance, 
tracking an infeasible ex ante optimal choice closely in Monte Carlo experiments.

We illustrate the framework using the cryptocurrency data described above. We estimate occupational characteristics of spot volatility that capture its average level, its variation, and its upper tail. The TRUMP token stands out on all three dimensions, combining a much higher average volatility level than BTC with greater variation and a far heavier upper tail. The occupational correlation between the log-price and spot volatility is clearly negative for BTC, in line with leverage-type dependence, whereas TRUMP displays a positive correlation, reflecting the joint post-launch decline of its price and volatility levels. Finally, an event analysis based on Federal Open Market Committee announcements and on Truth Social posts classified as financially relevant by a large language model (LLM) shows that spot volatility is elevated around information arrivals, with policy announcements acting market-wide and Truth Social posts operating mainly through a TRUMP-specific attention channel.

The rest of the paper is organized as follows. Section~\ref{sec:spot-vol} introduces the model and establishes a uniform approximation result for the spot volatility estimator. Section~\ref{sec:functionals} develops the uniform estimation theory for integrated volatility functionals and proposes a data-driven choice of the power index. Section~\ref{sec:simulation} reports the Monte Carlo evidence, and Section~\ref{sec:empirical} presents the empirical application. Section~\ref{sec:conclusion} concludes. All proofs, together with the details of the LLM-based event classification, are collected in the online supplementary materials.

\section{Spot volatility estimation in pure-jump models}
\label{sec:spot-vol}

In Section~\ref{subsec:settings}, we introduce the underlying pure-jump semimartingale model for the price process. In Section~\ref{subsec:uniform-spot}, we define the block-based spot volatility estimator and establish a uniform approximation result that holds across all blocks.

Throughout the paper, we use $\xrightarrow{\mathbb{P}}$ and $\xrightarrow{\mathbb{P}^*}$ to denote convergence in probability and convergence in outer probability, respectively.
For a sequence of random variables $\{X_n\}$, we write $X_n=o_p(1)$ if $X_n\xrightarrow{\mathbb{P}}0$.
For two real sequences $\{a_n\}$ and $\{b_n\}$, we write $a_n\asymp b_n$ if there exists a constant $K\ge 1$ such that
$a_n/K \le b_n \le K a_n$.
We also write 
$a\wedge b\equiv\min\{a,b\}$ and $(a)_+\equiv\max\{a,0\}$.
Let $\mathbb{R}$ be the set of real numbers,  $\mathbb{R}_+\equiv(0,\infty)$, and $\mathbb{R}^q$ the $q$-dimensional Euclidean space equipped with the Euclidean norm $\|\cdot\|$.
For $a\in\mathbb{R}$, $\lfloor a\rfloor$ denotes the greatest integer not exceeding $a$.
All limits are taken as $n\to\infty$.

\subsection{The settings}
\label{subsec:settings}
Let $X$ be a continuous-time pure-jump semimartingale price process defined on a filtered probability space \((\Omega, \mathcal{F}, (\mathcal{F}_t)_{t \geq 0}, \mathbb{P})\), given by
\begin{align}
\label{eq:price-process}
X_t = X_0 + \int_0^t b_s \,ds + \int_0^t \sigma_{s-}\,dZ_s,
\end{align}
where \(b\) is a drift process, 
 \(\sigma\) is a (positive) volatility  process with c\`adl\`ag sample path, and \(Z\) is a symmetric \(\beta\)-stable process with jump-activity index \(\beta \in (1,2)\), satisfying for any \(t \geq s\) and \(u \in \mathbb{R}\),
\begin{align*}
	\mathbb{E}\!\left[e^{iu(Z_t-Z_s)}\mid \mathcal{F}_s\right] 
	= \exp\!\left(-\tfrac{1}{2} (t-s)|u|^{\beta}\right).
\end{align*}
When $\beta=2$, $Z$ reduces to a standard Brownian motion and \eqref{eq:price-process} becomes the classical It\^{o} semimartingale model with a continuous diffusion component, which we regard as the standard benchmark. Our analysis complements it by developing the estimation theory for the heavy-tailed regime $\beta\in(1,2)$.

Assumption~\ref{assumption} collects mild regularity conditions on the drift and volatility processes.
\begin{assumption}
\label{assumption}
Suppose that there exists a sequence \((T_m)_{m \geq 1}\) of stopping times increasing to infinity and a sequence \((K_m)_{m \geq 1}\) of constants such that the following conditions hold for each \(m \geq 1\): (i) 
\(
|b_t| + |\sigma_t|+|\sigma_t|^{-1} \leq K_m\,\,  \text{for all } t \in [0, T_m]
\); (ii)
\(
\mathbb{E}[|\sigma_{t \wedge T_m} - \sigma_{s \wedge T_m}|^2] \leq K_m |t-s|^{2\kappa}\,\,  \text{for all } t, s \ge0
\) and for some $\kappa>0$.
\end{assumption}

Condition (i) requires the drift and volatility processes to be locally bounded, which is a standard assumption in high-frequency econometrics and statistics (see, for example, \citealp{jacod2012discretization,ait2014high}).
In addition, this condition implies that $\sigma_t$ takes values in a compact 
subset of $\mathbb R_+$ on each $[0, T_m]$.
Condition (ii) imposes local $\kappa$-H\"{o}lder continuity of the volatility process in the $L^2$ norm for some $\kappa>0$.
The index $\kappa$ may be arbitrarily small. The condition holds with
$\kappa=1/2$ when $\sigma$ is an It\^{o} semimartingale or a long-memory
process driven by fractional Brownian motion (\citealp{comte1996long}),
and smaller values of $\kappa$ accommodate rough volatility models
(\citealp{gatheral2018volatility,chong2024statistical}).

Assume that $X$ is sampled at the regularly spaced times $i\Delta_n$,
$i=0,1,\ldots,\lfloor T/\Delta_n\rfloor$, over a fixed interval $[0,T]$,
where the sampling interval $\Delta_n$ shrinks to zero in the asymptotic
scheme.
The $i$th increment of $X$ is
defined as 
\begin{align*}
\Delta_i^nX\equiv X_{i\Delta_n} - X_{(i-1)\Delta_n}.
\end{align*}
We choose a sequence of integers $k_n$ satisfying
$k_n\asymp\Delta_n^{-\gamma}$ for some $\gamma\in(0,1)$, and   
we divide the sample into $\lfloor T/(k_n\Delta_n)\rfloor$ non-overlapping 
blocks, each of which contains $k_n$ returns. 
Let \(\mathcal{I}_{n} \equiv \{0, \ldots, \lfloor T/(k_n\Delta_n)\rfloor-1\}\) denote the collection of block indices, with the $i$th block spanning the time interval $[ik_n\Delta_n,(i+1)k_n\Delta_n)$.

\subsection{A uniform approximation result for spot volatility estimators}
\label{subsec:uniform-spot}

For the $i$th block, the spot volatility estimator is defined as 
\begin{align}
\label{spot-estimator}
\widehat{\sigma}^p_{n,i}\equiv \dfrac{1}{c_{\beta}(p)k_n\Delta_n^{p/\beta}}\sum_{j=1}^{k_n} |\Delta_{ik_n+j}^nX|^p,\,\,0<p<\beta,
\end{align}
where \(c_{\beta}(p)\equiv\mathbb{E}[|Z_1|^p]\) is a normalizing constant, whose explicit form is (see page~163 of \citealp{ken1999levy})
\begin{align}
\label{eq:normal-const}
	c_{\beta}(p)
    =
    \frac{
        2^{p - p/\beta}
        \Gamma\!\left(\dfrac{1+p}{2}\right)
        \Gamma\!\left(1 - \dfrac{p}{\beta}\right)
    }{
        \sqrt{\pi}\,
        \Gamma\!\left(1 - \dfrac{p}{2}\right)
    }.
\end{align}
The estimator depends on the exponent $p$, which is held fixed at a value in $(0,\beta)$ throughout our theoretical development.\footnote{An alternative approach estimates volatility through
the empirical characteristic function of high-frequency increments,
with the frequency at which the characteristic function is evaluated
playing a role analogous to that of the power index $p$ here. This
approach has been developed for integrated volatility functionals,
including the locally stable pure-jump case
(\citealp{todorov2012realized,todorov2012realizeda}), and for spot
volatility in models with a Brownian component
(\citealp{JacodTodorov2014,JacodTodorov2018}), but not yet for spot
volatility estimation in the stable setting considered here.}
Accordingly, we suppress this dependence in the notation, and we discuss the practical choice of $p$ in Section~\ref{subsec:algo}. Let $\bar{\sigma}^p_{n,i}$ be the local average of
$|\sigma_t|^p$ over the window $[ik_n\Delta_n,(i+1)k_n\Delta_n)$, defined by
\begin{align}
\label{eq:moving-average}
	\bar{\sigma}^p_{n,i} \equiv \dfrac{1}{k_n\Delta_n}\int_{ik_n\Delta_n}^{(i+1)k_n\Delta_n} |\sigma_s|^p\,ds,\,\,0<p<\beta.
\end{align} 

The following theorem shows that $\widehat{\sigma}^p_{n,i}$ approximates the local average $\bar{\sigma}^p_{n,i}$, both for each fixed block and uniformly across all blocks.

\begin{theorem}
\label{th:spot-uniform}
For any fixed $p\in(0,\beta)$ and $k_n\asymp\Delta_n^{-\gamma}$:
\begin{enumerate}
\item[(i)] If Assumption~\ref{assumption}(i) holds and $\gamma\in(0,1)$, then for any 
fixed $i\in\mathcal{I}_n$,
\(
\bigl|\widehat{\sigma}^p_{n,i}-\bar{\sigma}^p_{n,i}\bigr|\xrightarrow{\mathbb{P}}0.
\)
\item[(ii)] If Assumption~\ref{assumption}(i) and  Assumption~\ref{assumption}(ii) hold and $\gamma$ satisfies
\begin{align}
\label{eq:gamma-range}
\max\left\{\frac{p}{\beta}+(1-p)_+,\;\;\frac{1}{2},\;\;1-(p\wedge 1)\kappa\right\}<\gamma<1,
\end{align}
then
\[
\sup_{i\in\mathcal{I}_n}\bigl|\widehat{\sigma}^p_{n,i}-\bar{\sigma}^p_{n,i}\bigr|
\xrightarrow{\mathbb{P}}0.
\]
\end{enumerate}
\end{theorem}

The first part of Theorem~\ref{th:spot-uniform} establishes the 
approximation for each fixed block under only the local boundedness 
condition in Assumption~\ref{assumption}(i), and the second part, which 
is the main result used in the sequel, establishes the approximation 
uniformly over all blocks under the additional H\"older regularity in 
Assumption~\ref{assumption}(ii). The stronger restriction on $\gamma$ 
in the second part ensures that the number of blocks, of order 
$T/(k_n\Delta_n)$, does not grow too quickly, and it provides sufficient 
control over the maximal estimation error across blocks. The lower bound 
on $\gamma$ in \eqref{eq:gamma-range} reflects three sources of 
estimation error: the drift contribution, the maximal fluctuation of the 
stable increments across blocks, and the local approximation of the spot 
volatility within each block.

We emphasize that Theorem~\ref{th:spot-uniform} concerns the 
approximation of the local averages $\bar{\sigma}^p_{n,i}$, rather than 
of the spot volatility process itself. Approximating the spot process 
uniformly would require the volatility path to be continuous within 
every estimation block, which, as the blocks shrink, essentially rules 
out volatility jumps. No such condition is needed here, as the $L^2$ 
H\"older condition in Assumption~\ref{assumption}(ii) is compatible with 
discontinuous volatility paths. This weaker notion of approximation 
suffices for our purposes, because the results for integrated volatility 
functionals in Section~\ref{sec:functionals} are established through 
a spatial localization technique that operates directly on the local 
averages.

\section{Uniform estimation of integrated volatility functionals}
\label{sec:functionals}

In this section, we define integrated volatility functionals through the occupation measure of the volatility process and construct their empirical counterparts by plugging in the spot volatility estimator from Section~\ref{sec:spot-vol}. Section~\ref{subsec:pointwise} establishes pointwise consistency for general random test functions, and Sections~\ref{subsec:bounded}--\ref{subsec:holder} establish uniform consistency over the bounded monotone, Lipschitz-in-parameter, and H\"older classes, with applications to volatility occupation times and quantiles and to argmax consistency of $M$-estimators. Section~\ref{subsec:algo} proposes a data-driven choice of the power index $p$.

Let $\mathbb{F}_T$ be the occupation measure induced by $\sigma_t$ over the time interval $[0,T]$, defined as (see, for example, \citealp{geman1980occupation,li2013volatility,li2021glivenko})
\begin{equation*}
\mathbb{F}_T(A) \equiv \int_0^T \mathbf1_{\{\sigma_s \in A\}}\, ds ,
\end{equation*}
for any Borel set $A \subseteq \mathbb{R}_+$. This measure records the amount of time that the spot volatility process spends in the set $A$.
By standard integration theory \citep{folland1999real}, $\mathbb{F}_T$ can equivalently be viewed as a linear functional acting on measurable functions $f:\mathbb{R}_+\to\mathbb{R}$, namely,
\begin{align}
\label{eq:inte-vol-funct}
   \mathbb{F}_{T}f\equiv \int_{\mathbb{R}_+} f(x)\,\mathbb{F}_T(dx)
   = \int_0^T f(\sigma_s)\,ds.
\end{align}
We refer to $\mathbb{F}_Tf$ as the integrated volatility functional associated 
with the test function $f$. When normalized by $T$, the occupation 
measure is a probability measure on $\mathbb{R}_+$, and $T^{-1}\mathbb{F}_Tf$ 
is the moment of $f$ under this measure, with averaging over time playing 
the role of the expectation. Classical distributional summaries, such as 
variances, covariances, and correlations, therefore have natural 
occupational counterparts, and we exploit this analogy in the 
applications below.

The empirical counterpart of $\mathbb{F}_T$ is constructed from the spot volatility estimator in \eqref{spot-estimator}.
To approximate the volatility path within each block, for $0\le i \le \lfloor T/(k_n\Delta_n)\rfloor-1$, define the piecewise-constant process
\begin{equation}
\label{eq:piecewise-process}
\begin{cases}
\widehat{\sigma}_{n,t} \equiv (\widehat{\sigma}^p_{n,i})^{1/p}, \quad t\in[ik_n\Delta_n,(i+1)k_n\Delta_n),\\
\widehat{\sigma}_{n,t} \equiv (\widehat{\sigma}^p_{n,\lfloor T/(k_n\Delta_n)\rfloor-1})^{1/p}, \quad t\in[\lfloor T/(k_n\Delta_n)\rfloor k_n\Delta_n,T].
\end{cases}
\end{equation}
The corresponding empirical occupation measure $\widehat{\mathbb{F}}_{n,T}$ is defined by
\begin{align}
\label{eq:inte-vol-funct-estimator}
\widehat{\mathbb{F}}_{n,T} f \equiv \int_0^T f(\widehat{\sigma}_{n,s})\,ds.
\end{align}

\subsection{Pointwise convergence for general test functions}
\label{subsec:pointwise}

In this subsection, we establish a general pointwise consistency result for  
$\widehat{\mathbb{F}}_{n,T}$, allowing the test function $f$ to be random. For this
purpose, we impose the following condition on $f$.
\begin{assumption}
\label{assumption-2}
(i) For Lebesgue almost every $t\in[0,T]$, $f$ is almost surely continuous at $\sigma_t$, i.e.,
    \(
        \mathbb{P}[\{\omega: f(\omega,\cdot)\ \text{is continuous at }\sigma_t(\omega)\}]=1.
    \)
     (ii) For a deterministic function $F:\mathbb{R}_+\to\mathbb{R}_+$, \(
     |f(\omega,x)|\le F(x)\,\text{for all }\omega\in\Omega, x\in\mathbb{R}_+.
    \)
    Moreover, $\sup_{x\in \mathcal{K}}F(x)<\infty$ for any compact subset $\mathcal{K}\subseteq\mathbb{R}_+$.
\end{assumption}

Condition (i) imposes a mild continuity requirement on the random function $f(\omega,\cdot)$, ensuring that any discontinuities of $f$ do not occur at values visited by the volatility path. Condition (ii) requires $f$ to admit a deterministic envelope $F$ that is bounded on compact sets. For uniformly bounded $f$, one may simply take $F$ to be a constant.

\begin{theorem}
\label{th:OM-point-consist}
 For any fixed $p\in (0,\beta)$,
 suppose that the following conditions hold: 
(i) Assumption \ref{assumption} and Assumption \ref{assumption-2}; 
(ii)  $k_n\asymp\Delta_n^{-\gamma}$ for some $\gamma$ satisfying \eqref{eq:gamma-range}.
Then 
    \begin{align*}
      \widehat{\mathbb{F}}_{n,T} f \xrightarrow{\mathbb{P}} \mathbb{F}_T f.
    \end{align*}
\end{theorem}

Pointwise convergence results of this kind are well established when the 
price process is an It\^{o} semimartingale with a Brownian diffusion 
component; see, for example, Theorem~9.4.1 of \citet{jacod2012discretization}, 
Lemma~1 of \citet{li2013volatility}, Theorem~3 of \citet{li2017adaptive}, and 
Theorem~2 of \citet{li2021glivenko}. Theorem~\ref{th:OM-point-consist} 
extends these results to the pure-jump setting. In applications, the 
test function $f$ is typically nonrandom. Allowing $f$ to be random is 
nevertheless a useful theoretical device, because the bracketing 
argument underlying our uniform convergence proofs in the subsequent 
sections involves random test functions as brackets, and 
Theorem~\ref{th:OM-point-consist} supplies the pointwise convergence 
required for these random brackets. With this building block in 
place, we now turn to the main results of this section, namely 
Glivenko--Cantelli type uniform convergence theorems for the empirical 
occupation measure over various classes of test functions.

\subsection{Uniform convergence for the class of bounded monotone functions}
\label{subsec:bounded}

Let $\mathcal{G}_M\equiv\{g:\mathbb{R}_+\to[0,1]\,:\,g\text{ is monotone}\}$ 
denote the class of bounded monotone functions on $\mathbb{R}_+$, and 
$\mathcal{G}_I\equiv\{\mathbf1_{(0,x]}:x\in\mathbb{R}_+\}\subseteq\mathcal{G}_M$ 
the subclass of indicator functions. The latter is of particular
interest: for each $x\in\mathbb{R}_+$, the volatility occupation time
\begin{equation}
\label{eq:occupation-time}
F_T(x)\equiv\mathbb{F}_T \mathbf1_{(0,x]} = \int_0^T \mathbf1_{\{\sigma_s\le x\}}\,ds
\end{equation}
records the amount of time over $[0,T]$ that $\sigma_t$ spends at or below 
level $x$. The normalized occupation time $F_T(x)/T$ is thus the pathwise 
analogue of a cumulative distribution function, and uniform estimation 
over $\mathcal{G}_I$ is the occupation-measure analogue of the classical 
Glivenko--Cantelli problem.

\begin{theorem}
\label{th:OM-uniform-consist-G_M}
 For any fixed $p\in (0,\beta)$,
     suppose that the following conditions hold: (i) Assumption \ref{assumption}; 
     (ii) $k_n\asymp\Delta_n^{-\gamma}$ for some $\gamma$ 
    satisfying \eqref{eq:gamma-range};
     (iii) the map $x\mapsto F_T(x)$ is almost surely continuous. Then  
    \begin{align*}
      \sup_{g\in\mathcal{G}_M} |\widehat{\mathbb{F}}_{n,T} g - \mathbb{F}_{T} g|\xrightarrow{\mathbb{P}^*}0.
    \end{align*}
\end{theorem}

The use of outer probability $\mathbb{P}^*$ avoids measurability issues 
that may arise when taking a supremum over an uncountable class (see 
\citealp{vdvw1996weak}). The proof relies on a random bracketing 
argument, and condition~(iii) ensures that the pointwise convergence 
of Theorem~\ref{th:OM-point-consist} applies to the random bracket 
functions constructed in this argument. 
Theorem~\ref{th:OM-uniform-consist-G_M} has two immediate consequences 
for the estimation of volatility occupation times and their quantiles, 
which we develop next.

Specializing \eqref{eq:inte-vol-funct-estimator} to $f=\mathbf1_{(0,x]}$ yields
the empirical volatility occupation time
\begin{align}
\label{eq:empirical-occupation-time}
\widehat{F}_{n,T}(x)\equiv\int_0^T \mathbf1_{\{\widehat{\sigma}_{n,s}\le x\}}\,ds.
\end{align}
Corollary~\ref{cor:occupation-time} establishes its uniform consistency over $x\in\mathbb{R}_+$.
\begin{corollary}
\label{cor:occupation-time}
Suppose that the conditions of Theorem~\ref{th:OM-uniform-consist-G_M} hold. Then
\begin{align*}
\sup_{x\in\mathbb{R}_+} \big|\widehat{F}_{n,T}(x) - F_T(x)\big|\xrightarrow{\mathbb{P}}0.
\end{align*}
\end{corollary}

By analogy with the classical definition of quantiles for probability 
distributions, we define the volatility occupation quantile as the 
left-continuous generalized inverse of $F_T(\cdot)$ (see, for example, 
\citealp{li2013volatility}). Specifically, for any $\alpha\in(0,1)$,
\begin{align}
\label{volatility-quantile}
Q_T(\alpha) \equiv \inf\{x\in\mathbb{R}_+:\,F_T(x)\geq \alpha T\}.
\end{align}
Its empirical counterpart is defined as
\begin{align}
\label{volatility-quantile-estimator}
\widehat{Q}_{n,T}(\alpha) \equiv \inf\{x\in\mathbb{R}_+:\,\widehat{F}_{n,T}(x)\geq \alpha T\}.
\end{align}
Corollary~\ref{cor:volatility-quantile} establishes its consistency at each continuity point of $Q_T(\cdot)$. This restriction is mild and standard, as it only rules out the flat regions of $F_T$ at which the quantile is not single-valued.
\begin{corollary}
\label{cor:volatility-quantile}
Let
\(
\mathcal{Q}\equiv\{\alpha \in (0,1):Q_T\ \text{is almost surely continuous at}\ \alpha\}.
\)
Suppose that the conditions of Theorem~\ref{th:OM-uniform-consist-G_M} hold. Then, for any fixed $\alpha \in \mathcal{Q}$,
\begin{align*}
\widehat{Q}_{n,T}(\alpha) \xrightarrow{\mathbb{P}} Q_T(\alpha).
\end{align*}
\end{corollary}

\subsection{Uniform convergence for the class of Lipschitz-in-parameter functions}
\label{subsec:lipschitz}

We now turn to classes of test functions indexed by a finite-dimensional 
parameter, which arise naturally in the estimation of parametric models. 
For such classes, it is particularly useful to let the test function 
depend also on an observable state process, since parametric criteria 
typically involve observed variables alongside the latent volatility. 
Let $Y\equiv\{Y_t:t\in[0,T]\}$ be an $\mathbb{R}^q$-valued
c\`adl\`ag process defined on $(\Omega, \mathcal{F}, (\mathcal{F}_t)_{t\ge 0}, \mathbb{P})$,
observed at the same sampling times $i\Delta_n$ as $X$, so that the data
consist of the discrete record
$\{Y_{i\Delta_n}:i=0,1,\ldots,\lfloor T/\Delta_n\rfloor\}$.
The process $Y$ may collect state variables relevant for time-varying 
market conditions, derivative pricing, or volatility patterns; see, for example, 
\citet{andersen1996return,avramov2006impact,chordia2001market,andersen2001variance}.

Specifically, let $\Theta\subseteq\mathbb{R}^{d}$ be a compact parameter space and 
$g:\mathbb{R}_+\times\mathbb{R}^q\times\Theta\to\mathbb{R}$ a measurable map. 
For each $\theta\in\Theta$ and $(x,y)\in\mathbb{R}_+\times\mathbb{R}^{q}$, set $g_\theta(x,y)\equiv g(x,y;\theta)$, and 
denote $\mathcal{G}_L\equiv\{g_\theta:\theta\in\Theta\}$. We assume that 
$(x,y)\mapsto g_\theta(x,y)$ is continuous for each $\theta\in\Theta$, 
and that the family is locally Lipschitz in $\theta$: for every compact 
$\mathcal{K}\subseteq\mathbb{R}_+\times\mathbb{R}^q$, there exists 
$M_\mathcal{K}>0$ such that
\[
\sup_{(x,y)\in\mathcal{K}}\bigl|g_\theta(x,y)-g_{\theta'}(x,y)\bigr|
\le M_\mathcal{K}\|\theta-\theta'\|, \quad \theta,\theta'\in\Theta.
\]

The empirical functional is constructed from these discrete
observations. In analogy with \eqref{eq:piecewise-process}, let
$Y_{n,t}\equiv Y_{ik_n\Delta_n}$ for $t\in[ik_n\Delta_n,(i+1)k_n\Delta_n)$
and $i\in\mathcal{I}_n$, with the value on the terminal segment
$[\lfloor T/(k_n\Delta_n)\rfloor k_n\Delta_n,T]$ given by that of the
last block. Extending the notation in \eqref{eq:inte-vol-funct}, we
define, for any $g_\theta\in\mathcal{G}_L$,
\begin{align*}
{\mathbb{F}}_{T} g_\theta \equiv \int_0^T g_\theta({\sigma}_{s},Y_s)\,ds, 	\qquad
\widehat{\mathbb{F}}_{n,T} g_\theta \equiv \int_0^T g_\theta(\widehat{\sigma}_{n,s},Y_{n,s})\,ds.
\end{align*}

\begin{theorem}
\label{th:OM-uniform-consist-G_L}
For any fixed $p\in(0,\beta)$, suppose that the following conditions hold:
(i) Assumption~\ref{assumption};
(ii) there exists a sequence $(\mathcal{K}'_m)_{m\ge1}$ of compact subsets of $\mathbb{R}^q$ such that $Y_t\in\mathcal{K}'_m$ for all $t\leq T_m$, where $(T_m)_{m\ge1}$ are the stopping times in Assumption~\ref{assumption};
(iii) $k_n\asymp\Delta_n^{-\gamma}$ for some $\gamma$ satisfying \eqref{eq:gamma-range}. Then
\begin{align*}
\sup_{g_\theta\in\mathcal{G}_L} |\widehat{\mathbb{F}}_{n,T} g_\theta - \mathbb{F}_T g_\theta|\xrightarrow{\mathbb{P}^*}0.
\end{align*}
\end{theorem}

Theorem~\ref{th:OM-uniform-consist-G_L} establishes uniform consistency 
of $\widehat{\mathbb{F}}_{n,T}$ over $\mathcal{G}_L$. In contrast to 
Theorem~\ref{th:OM-uniform-consist-G_M}, no continuity assumption on the 
occupation time function is required, because the test functions in 
$\mathcal{G}_L$ are themselves continuous and the random bracketing argument 
can be replaced by a deterministic finite-cover argument over $\Theta$.

Theorem~\ref{th:OM-uniform-consist-G_L} can be used to establish consistency 
of $M$-estimators within the standard argmax consistency framework 
(see, for example, Chapter~5 of \citealp{vdv1998asymptotic}). In particular, it provides uniform convergence of the sample criterion over the parameter space. Let $\theta\in\Theta$ denote the finite-dimensional parameter of interest, and define the sample criterion 
function
\begin{equation}\label{eq:criterion-function}
M_n(\theta)\;\equiv\;\widehat{\mathbb F}_{n,T} m_\theta
=\int_0^T m_\theta(\widehat\sigma_{n,s},Y_{n,s})\,ds,
\end{equation}
where the collection of test functions $\{m_\theta:\theta\in\Theta\}$ is assumed 
to be contained in $\mathcal G_L$. The corresponding population criterion is
\[
M(\theta)\equiv \mathbb F_T m_\theta = \int_0^T m_\theta(\sigma_s,Y_s)\,ds, \quad \theta\in\Theta.
\]

Let $\theta_0$ denote the unique maximizer of $M$, and let $\widehat\theta_n$ 
be a sequence of (approximate) maximizers of $M_n$. 
Theorem~\ref{th:M-estimators} below establishes the consistency of 
$\widehat\theta_n$ for $\theta_0$.

\begin{theorem}
\label{th:M-estimators}
For any fixed $p\in(0,\beta)$,
suppose that the following conditions hold:
(i) the conditions of Theorem~\ref{th:OM-uniform-consist-G_L} hold;
(ii) for every $\epsilon>0$, 
$
\sup_{\theta\in\Theta:\|\theta-\theta_0\|\geq \epsilon} M(\theta)<M(\theta_0)
$ almost surely.
Then, for any sequence of estimators $\widehat{\theta}_n$ satisfying
$
M_n(\widehat{\theta}_n)\geq M_n(\theta_0)-o_p(1),
$
\begin{align*}
\widehat{\theta}_n\xrightarrow{\mathbb{P}^*}\theta_0.
\end{align*}
\end{theorem}

Condition~(i) ensures that the sample criterion converges to its population 
counterpart uniformly over $\Theta$ via Theorem~\ref{th:OM-uniform-consist-G_L}. 
Condition~(ii) is the standard identification requirement that $M$ has a unique 
and well-separated maximizer at $\theta_0$ (see, for example, Theorem~5.7 of 
\citealp{vdv1998asymptotic}). Compared to the classical $M$-estimation framework 
based on the i.i.d. empirical measure \citep{vdv1998asymptotic}, our setting 
replaces $P_n$ with the empirical occupation measure 
$\widehat{\mathbb{F}}_{n,T}$ constructed from the spot volatility estimator. 
We also note that this framework does not require the parametric model 
underlying $m_\theta$ to be correctly specified. Under misspecification, 
$\theta_0$ is interpreted as a pseudo-true parameter in the spirit of 
\citet{white1982maximum}, defined pathwise as the maximizer of the 
occupational criterion $M$.

\begin{example}[Derivative pricing]
\label{ex:derivative-pricing}
\citet{li2016generalized} study the estimation of derivative pricing models using integrated moment conditions when the underlying price has a Brownian component, and our theory extends this type of consistency result to the pure-jump setting. Specifically, let $P_t$ denote the observed price at time $t$ of a derivative written on an underlying asset whose time-$t$ price and latent spot volatility are $X_t$ and $\sigma_t$, respectively. Set $Y_t\equiv(T-t,X_t,P_t)$, an observable vector-valued process. Suppose that the derivative price satisfies $P_t=h(T-t,X_t,\sigma_t;\theta_0)$, where $h$ is a known pricing function implied by a parametric risk-neutral model (e.g., a stochastic volatility model) indexed by $\theta\in\Theta$, and $\theta_0\in\Theta$ is the true parameter. A natural criterion function $m_\theta$ is then given by $m_\theta(\sigma_t,Y_t)\equiv-\rho\bigl(P_t-h(T-t,X_t,\sigma_t;\theta)\bigr)$, where $\rho$ is a given loss function. If $\{m_\theta:\theta\in\Theta\}\subset\mathcal G_L$ and the identification condition (Theorem~\ref{th:M-estimators}(ii)) holds, then $\theta_0$ can be consistently estimated by maximizing the empirical criterion $M_n(\theta)$ defined in \eqref{eq:criterion-function}. If the pricing model is misspecified, $\theta_0$ is accordingly interpreted as the pseudo-true parameter associated with the loss $\rho$.
\end{example}

\begin{example}[Volatility pattern approximation]
\label{exp:volatility-pattern}
Intraday volatility is known to exhibit systematic time-of-day patterns; see, for example, \citet{andersen2001variance}. Let $Y_t\equiv t-\lfloor t\rfloor$ denote the time of day at time $t$, with one day as the unit of time, and let $\widetilde h(\cdot;\theta)$ be a parametric family of intraday pattern functions indexed by $\theta\in\Theta$. Consider the quadratic loss, under which $m_\theta(x,y)\equiv-(x-\widetilde h(y;\theta))^2$ and the population criterion becomes $M(\theta)=-\int_0^T (\sigma_s-\widetilde h(Y_s;\theta))^2\,ds$. The maximizer $\theta_0$ thus identifies the best approximation of the volatility path within the parametric family in the occupational $L^2$ sense, and the family need not be correctly specified. If $\{m_\theta:\theta\in\Theta\}\subset\mathcal G_L$ and the identification condition (Theorem~\ref{th:M-estimators}(ii)) holds, the estimated intraday pattern $\widetilde h(\cdot;\widehat\theta_n)$ consistently recovers this occupational projection.
\end{example}

Theorem~\ref{th:OM-uniform-consist-G_L} also delivers a useful 
by-product that relies only on pointwise convergence. Since uniform 
convergence implies pointwise convergence, the empirical occupation 
measure consistently estimates $\mathbb{F}_T g$ for each fixed 
$g\in\mathcal{G}_L$, and quantities that involve only finitely many such 
functionals therefore require no uniformity. The following example 
illustrates this point.

\begin{example}[Occupational correlation and leverage effect]
\label{exp:correlation}
A leading case is the leverage effect, commonly understood as a negative dependence between asset returns and volatility; see, for example, \citet{ait2017estimation,bandi2012time,chen2024leverage,kalnina2017nonparametric,wang2014estimation}. We take the state variable to be the price itself, that is, $Y_t\equiv X_t$. With $\bar{\sigma}\equiv T^{-1}\int_0^T \sigma_s\,ds$ and $\bar{Y}\equiv T^{-1}\int_0^T Y_s\,ds$, the occupational covariance and correlation between $\sigma$ and $Y$ are
\begin{align*}
\Cov_T(\sigma,Y)
\equiv
\int_0^T (\sigma_s-\bar{\sigma})(Y_s-\bar{Y})\,ds,
\qquad
\Corr_T(\sigma,Y)
\equiv
\frac{\Cov_T(\sigma,Y)}{\sqrt{\Var_T(\sigma)\,\Var_T(Y)}},
\end{align*}
where $\Var_T(\sigma)\equiv\int_0^T(\sigma_s-\bar{\sigma})^2\,ds$ and $\Var_T(Y)\equiv\int_0^T(Y_s-\bar{Y})^2\,ds$. All ingredients involving $\sigma$ are integrated functionals generated by the three test functions $g(x,y)\in\{x,\,xy,\,x^2\}$, so pointwise convergence suffices. Combined with the continuous mapping theorem, it delivers consistent plug-in estimation of $\Corr_T(\sigma,Y)$. This occupational correlation provides a global measure of the dependence between volatility and prices, in the spirit of the early leverage-effect studies that related volatility to observed price levels \citep{black1976studies,christie1982stochastic}. It is distinct from the instantaneous volatility--price covariations studied in the existing high-frequency literature, which are local in nature and are tied to the Brownian shocks of the underlying processes. The two are complementary, capturing the same co-movement between volatility and prices on different time scales.
\end{example}

\subsection{Uniform convergence for 
the class of H\"{o}lder continuous functions}
\label{subsec:holder}

Let $\mathcal{G}_H$ be a class of real-valued functions on $\mathbb{R}_+$ satisfying the following condition: for any compact subset $\mathcal{K}\subseteq\mathbb{R}_+$, there exist constants $\alpha_{\mathcal K}\in(0,1]$ and $M_{\mathcal K}>0$ such that
\[
\sup_{x,y\in\mathcal K,\,x\ne y}
\frac{|g(x)-g(y)|}{|x-y|^{\alpha_{\mathcal K}}}
\le M_{\mathcal K},
\qquad g\in\mathcal G_H.
\]
That is, the functions in $\mathcal{G}_H$ are uniformly $\alpha_{\mathcal K}$-H\"older continuous on each $\mathcal K$. 

\begin{theorem}
\label{th:OM-uniform-consist-G_H}
For any fixed $p\in(0,\beta)$, suppose that the following conditions hold:
(i) Assumption~\ref{assumption};
(ii) $k_n\asymp\Delta_n^{-\gamma}$ for some $\gamma$ satisfying \eqref{eq:gamma-range}. Then
\begin{align*}
\sup_{g\in\mathcal{G}_H} \big|\widehat{\mathbb{F}}_{n,T} g - \mathbb{F}_T g\big|\xrightarrow{\mathbb{P}^*}0.
\end{align*}
\end{theorem}

Theorem~\ref{th:OM-uniform-consist-G_H} establishes uniform consistency of 
the empirical occupation measure over the class of locally H\"{o}lder continuous 
test functions. This covers a broad family of nonlinear transformations commonly used in 
the study of integrated volatility functionals, including power functions 
(e.g., integrated variance or quarticity, \citealp{jacod2013quarticity}) and
exponential functions (e.g., volatility Laplace transforms,
\citealp{todorov2012realized,todorov2012realizeda}).

\begin{example}[Occupational moments]
For any $r>0$, let $g(x)\equiv x^r$. Since $g\in\mathcal G_H$, 
Theorem~\ref{th:OM-uniform-consist-G_H} yields 
$\widehat{\mathbb F}_{n,T} g\xrightarrow{\mathbb P}\mathbb F_T g$. 
Hence, the empirical occupation measure permits consistent estimation of 
occupational moments of arbitrary order. By the continuous mapping 
theorem, classical distributional summaries constructed from these 
moments, such as the occupational mean, variance, skewness, and kurtosis, can 
be consistently estimated whenever the corresponding population quantities 
are well defined.
\end{example}

\subsection{A data-driven choice of the power index}
\label{subsec:algo}

The pure-jump setting requires a choice of the power index $p$,
and the familiar choice $p=2$ is unavailable because the increments
have infinite variance. We recommend choosing $p$ based on a stability
consideration. Since the target functional does not depend on $p$, the 
estimate should be insensitive to small changes in $p$ within the 
range where consistency holds, and we search for a value of $p$ at 
which the estimate is stable in this sense.

The search procedure, summarized in the algorithm below, takes an 
initial lower bound $p_1$, a step size $\delta$, and a relative tolerance 
$\tau$ as inputs. The algorithm starts from a wide bracket of 
candidate values and compares the estimates computed at its two 
endpoints. If the two estimates agree within the tolerance $\tau$, the 
search stops and the midpoint of the bracket is returned. Otherwise, 
the bracket is shrunk by one step from the end whose removal brings 
the endpoint estimates closer together, and the comparison is 
repeated. The search is confined to $p<\beta/2$, where the spot 
volatility estimator attains the optimal convergence rate (see, for 
example, \citealp{yan2026nonparametric}), and it depends on $\beta$ only 
through the initialization of $p_2$, so that a preliminary estimate of 
$\beta$ suffices in applications. In practice, this stability-based 
choice delivers estimation results that are robust to the choice of $p$, 
and the simulation study in Section~\ref{sec:simulation} further 
suggests that it is effective in reducing estimation error.

\medskip
\noindent\fbox{\parbox{0.96\textwidth}{
\textbf{Algorithm} (Stability-based $p$-selection). \\[4pt]
\textbf{Input:} Increments $\Delta_i^n X$; activity index $\beta$ (a preliminary estimate suffices); test function
$g$; initial lower bound $p_1$; step size $\delta$; tolerance $\tau$. Here, 
$\widehat{\mathbb F}_{n,T}^{(p)}g$ denotes the estimator computed with power index $p$.\\[3pt]
\textbf{Initialize:}
$p_2\leftarrow \lfloor\beta/(2\delta)\rfloor\cdot\delta$,\;
$\widehat{F}_1\leftarrow \widehat{\mathbb F}_{n,T}^{(p_1)}g$,\;
$\widehat{F}_2\leftarrow \widehat{\mathbb F}_{n,T}^{(p_2)}g$.\\[3pt]
\textbf{While} $p_1+\delta<p_2$\textbf{:}
\begin{enumerate}
\item[1.] If $|\widehat{F}_1-\widehat{F}_2|/\max(|\widehat{F}_1|,|\widehat{F}_2|)<\tau$, \textbf{return} $(p_1+p_2)/2$.
\item[2.] Compute $p_1^+=p_1+\delta$, $p_2^-=p_2-\delta$.
\item[3.] Evaluate $\widehat{F}_1^+=\widehat{\mathbb F}_{n,T}^{(p_1^+)}g$, $\widehat{F}_2^-=\widehat{\mathbb F}_{n,T}^{(p_2^-)}g$.
\item[4.] Compute relative differences $r_A=|\widehat{F}_1^+-\widehat{F}_2|/\max(|\widehat{F}_1^+|,|\widehat{F}_2|)$, $r_B=|\widehat{F}_1-\widehat{F}_2^-|/\max(|\widehat{F}_1|,|\widehat{F}_2^-|)$.
\item[5.] If $r_A\le r_B$, set $p_1\leftarrow p_1^+$, $\widehat{F}_1\leftarrow \widehat{F}_1^+$; otherwise, set $p_2\leftarrow p_2^-$, $\widehat{F}_2\leftarrow \widehat{F}_2^-$.
\end{enumerate}
\textbf{End while.}\\
\textbf{Return} $(p_1+p_2)/2$.
}}
\medskip

\section{Simulation study}
\label{sec:simulation}

This section examines the finite-sample performance of the proposed 
estimators. Section~\ref{subsec:design} describes the Monte Carlo 
design, and Section~\ref{subsec:sim-results} reports the 
results.

\subsection{Monte Carlo design}
\label{subsec:design}

We generate data from the stochastic volatility model
\begin{equation}
\label{eq:SV-1}
dX_t = \sqrt{V_t}\, dZ_t,
\qquad
dV_t = 0.03\,(1 - V_t)\,dt + 0.2\sqrt{V_t}\, dB_t,
\end{equation}
where $Z$ is a symmetric $\beta$-stable process and $B$ is a standard 
Brownian motion independent of $Z$. The variance process $V$ follows a 
square-root diffusion calibrated as in \citet{li2013volatility}. The spot 
volatility is $\sigma_t=\sqrt{V_t}$. We initialize $V_0$ at the median of 
its invariant distribution and set the drift in $X$ to zero for simplicity.

We fix the time span at $T=22$ trading days (with one day as the unit of 
time), the intraday sample size at $n=390$ (corresponding to the 1-minute
sampling frequency over a 6.5-hour trading day), and the block size at 
$k_n=30$. The activity index is varied over 
$\beta\in\{1.1,1.2,\ldots,1.9\}$. We consider four continuous test functions 
$g(x)\in\{x^2, x, \sqrt{x}, \log x\}$ and three indicator test functions 
$g(x)=\mathbf 1_{\{x\le u\}}$ for $u\in\{u_{25},u_{50},u_{75}\}$, the empirical 
25\%, 50\%, and 75\% quantiles of $\{\sigma_t\}_{t\in[0,T]}$ in each 
replication. All results are based on 10{,}000 Monte Carlo replications.

For each $(\beta, g)$ pair, we evaluate the estimator over a grid of 
candidate values of $p$ and measure estimation accuracy by the mean 
absolute deviation (MAD) across replications, normalized by the absolute 
value of the true target and expressed in percent. We refer to this 
measure as the relative MAD. The true value of the functional is 
computed in each replication by numerically integrating the test 
function along the simulated volatility path.

We compare two choices of the power index. The first is the 
grid-optimal value $p^*$ that minimizes the relative MAD over the grid 
$\{0.01,0.02,\ldots,\beta\}$. Since computing $p^*$ requires the true 
value of the target functional, it is infeasible in practice and 
{serves as the ex ante optimal choice}.
The second is the data-driven choice 
$\widehat p^*_{\mathrm{algo}}$ obtained from the stability-based 
selection rule of Section~\ref{subsec:algo}, implemented with $p_1=0.05$, 
$\delta=0.01$, $\tau=0.1\%$, and the true value of $\beta$.

For each $(\beta,g)$ pair, we report the ex ante optimum $p^*$ and its 
relative MAD, together with the median of $\widehat p^*_{\mathrm{algo}}$ 
across replications and the relative MAD attained when each replication 
uses its own selected value. We summarize the comparison by the relative 
efficiency $\mathrm{MAD}(p^*)/\mathrm{MAD}_{\mathrm{algo}}$, with values 
close to one indicating that the selection rule performs comparably to the
{ex ante optimum}.

\subsection{Simulation results}
 \label{subsec:sim-results}

Table~\ref{tab:mc_continuous} reports the results 
for the four continuous test functions. For all four functions,  $p^*$ increases with $\beta$ and the associated relative 
MAD decreases with $\beta$. The convex case $g(x)=x^2$ is the most 
challenging, the square-root case yields the smallest errors, and the 
linear and logarithmic cases lie in between. The stability-based rule 
tracks 
the ex ante optimum 
closely across all functions and all levels of 
$\beta$, with  the resulting relative 
MADs nearly indistinguishable from their grid-optimal counterparts.

Table~\ref{tab:mc_indicator} reports the results for the indicator test 
functions. Finite-sample performance now depends more strongly on 
the location of the threshold. The lower-threshold case 
$g(x)=\mathbf 1_{\{x\le u_{25}\}}$ is the most difficult, with relative 
MAD above 26\% throughout, whereas the upper-threshold case attains the 
smallest MAD of roughly 4\%, and the median-threshold case lies 
in between. For all three thresholds, $p^*$ increases with $\beta$. The 
stability-based rule remains competitive for $u_{25}$ and $u_{50}$, and 
its gap relative to the {ex ante optimum} narrows as $\beta$ increases. The main 
exception is $u_{75}$, for which the rule searches over $p<\beta/2$ 
while the ex ante optimum  selects values above 
$\beta/2$, producing a more visible gap. Even in this case, the relative 
MAD of the data-driven choice remains moderate in absolute terms.

Table~\ref{tab:efficiency_all} summarizes the relative efficiency for 
all test functions. The stability-based rule comes close to, and in some 
cases slightly exceeds, the {ex ante optimum} for the smooth functions and for the 
indicator functions at thresholds $u_{25}$ and $u_{50}$. The 
high-threshold case remains the main exception, as discussed above.

Overall, the proposed estimators, combined with the data-driven 
choice of $p$, perform well for the commonly used test functions, with 
the stability-based rule closely tracking the {ex ante optimum} across all levels of $\beta$. The occupation time at a high 
threshold is an exception, although even in that case the estimation 
error remains moderate. We therefore recommend the procedure for 
practical use.

\begin{table}[htbp]
\centering
\small
\begin{threeparttable}
\caption{Monte Carlo results for continuous test functions}
\label{tab:mc_continuous}
\setlength{\tabcolsep}{4.5pt}
\begin{tabular}{c c c c c c c c c}
\toprule
 & \multicolumn{4}{c}{$g(x)=x^2$} & \multicolumn{4}{c}{$g(x)=x$} \\
\cmidrule(lr){2-5} \cmidrule(lr){6-9}
$\beta$ & $p^*$ & MAD (\%) & $\hat p^*_{\mathrm{algo}}$ & $\mathrm{MAD}_{\mathrm{algo}}$ (\%) & $p^*$ & MAD (\%) & $\hat p^*_{\mathrm{algo}}$ & $\mathrm{MAD}_{\mathrm{algo}}$ (\%) \\
\midrule
1.1 & 0.08 & 15.06 & 0.10 & 14.80 & 0.16 & 3.31 & 0.18 & 3.18 \\
1.2 & 0.12 & 13.10 & 0.14 & 12.77 & 0.21 & 2.86 & 0.23 & 2.74 \\
1.3 & 0.17 & 11.19 & 0.18 & 10.83 & 0.27 & 2.39 & 0.28 & 2.29 \\
1.4 & 0.22 & 9.71 & 0.24 & 9.32 & 0.33 & 2.03 & 0.33 & 1.95 \\
1.5 & 0.28 & 8.40 & 0.30 & 7.98 & 0.39 & 1.75 & 0.36 & 1.69 \\
1.6 & 0.33 & 7.24 & 0.38 & 6.77 & 0.46 & 1.49 & 0.42 & 1.46 \\
1.7 & 0.40 & 6.17 & 0.47 & 5.66 & 0.53 & 1.27 & 0.46 & 1.24 \\
1.8 & 0.49 & 5.18 & 0.59 & 4.66 & 0.61 & 1.09 & 0.52 & 1.08 \\
1.9 & 0.59 & 4.13 & 0.72 & 3.54 & 0.72 & 0.90 & 0.62 & 0.88 \\
\midrule
 & \multicolumn{4}{c}{$g(x)=\sqrt{x}$} & \multicolumn{4}{c}{$g(x)=\log(x)$} \\
\cmidrule(lr){2-5} \cmidrule(lr){6-9}
$\beta$ & $p^*$ & MAD (\%) & $\hat p^*_{\mathrm{algo}}$ & $\mathrm{MAD}_{\mathrm{algo}}$ (\%) & $p^*$ & MAD (\%) & $\hat p^*_{\mathrm{algo}}$ & $\mathrm{MAD}_{\mathrm{algo}}$ (\%) \\
\midrule
1.1 & 0.23 & 0.84 & 0.18 & 0.88 & 0.04 & 5.19 & 0.06 & 6.11 \\
1.2 & 0.27 & 0.75 & 0.21 & 0.78 & 0.05 & 4.90 & 0.08 & 5.55 \\
1.3 & 0.29 & 0.66 & 0.23 & 0.68 & 0.06 & 4.65 & 0.11 & 5.14 \\
1.4 & 0.32 & 0.59 & 0.26 & 0.60 & 0.08 & 4.41 & 0.14 & 4.77 \\
1.5 & 0.35 & 0.55 & 0.29 & 0.55 & 0.09 & 4.29 & 0.17 & 4.54 \\
1.6 & 0.39 & 0.50 & 0.33 & 0.50 & 0.12 & 4.13 & 0.21 & 4.29 \\
1.7 & 0.41 & 0.47 & 0.38 & 0.46 & 0.14 & 3.98 & 0.26 & 4.11 \\
1.8 & 0.47 & 0.43 & 0.43 & 0.42 & 0.19 & 3.85 & 0.32 & 3.93 \\
1.9 & 0.52 & 0.40 & 0.53 & 0.38 & 0.22 & 3.66 & 0.41 & 3.79 \\
\bottomrule
\end{tabular}
\begin{tablenotes}[flushleft]
\small
\item[] \textit{Note:} For each activity index $\beta$, $p^*$ is the ex ante optimum, i.e., the value of the power index that minimizes the relative MAD, computed using the true value of the target functional, and $\hat p^*_{\mathrm{algo}}$ is the median value selected by the stability-based rule across replications, with $\mathrm{MAD}_{\mathrm{algo}}$ the associated relative MAD. The relative MAD is the mean absolute deviation across replications, normalized by the absolute value of the true target and expressed in percent. Based on 10,000 Monte Carlo replications.
\end{tablenotes}
\end{threeparttable}
\end{table}

\begin{table}[htbp]
\centering
\footnotesize
\begin{threeparttable}
\caption{Monte Carlo results for indicator test functions}
\label{tab:mc_indicator}
\setlength{\tabcolsep}{1.4pt}
\renewcommand{\arraystretch}{1.0}
\begin{tabular}{c c c c c c c c c c c c c}
\toprule
& \multicolumn{4}{c}{$g(x)=\mathbf1_{\{x\le u_{25}\}}$}
& \multicolumn{4}{c}{$g(x)=\mathbf1_{\{x\le u_{50}\}}$}
& \multicolumn{4}{c}{$g(x)=\mathbf1_{\{x\le u_{75}\}}$} \\
\cmidrule(lr){2-5} \cmidrule(lr){6-9} \cmidrule(lr){10-13}
$\beta$ & $p^*$ & MAD (\%) & $\hat p^*_{\mathrm{algo}}$ & $\mathrm{MAD}_{\mathrm{algo}}$ (\%)
& $p^*$ & MAD (\%) & $\hat p^*_{\mathrm{algo}}$ & $\mathrm{MAD}_{\mathrm{algo}}$ (\%)
& $p^*$ & MAD (\%) & $\hat p^*_{\mathrm{algo}}$ & $\mathrm{MAD}_{\mathrm{algo}}$ (\%) \\
\midrule
1.1 & 0.02 & 40.25 & 0.15 & 44.80 & 0.03 & 7.24 & 0.17 & 8.49 & 0.67 & 3.77 & 0.22 & 11.14 \\
1.2 & 0.03 & 37.85 & 0.17 & 41.33 & 0.03 & 7.08 & 0.18 & 8.21 & 0.73 & 3.86 & 0.23 & 10.68 \\
1.3 & 0.04 & 36.23 & 0.18 & 38.74 & 0.03 & 7.11 & 0.20 & 8.11 & 0.78 & 3.91 & 0.24 & 10.13 \\
1.4 & 0.05 & 34.57 & 0.20 & 36.37 & 0.04 & 7.09 & 0.21 & 7.94 & 0.83 & 3.96 & 0.26 & 9.76 \\
1.5 & 0.07 & 33.39 & 0.23 & 34.39 & 0.05 & 7.10 & 0.23 & 7.84 & 0.88 & 3.96 & 0.28 & 9.32 \\
1.6 & 0.10 & 31.90 & 0.26 & 32.19 & 0.06 & 7.06 & 0.25 & 7.64 & 0.94 & 4.00 & 0.28 & 9.00 \\
1.7 & 0.14 & 30.68 & 0.30 & 30.64 & 0.07 & 7.00 & 0.28 & 7.47 & 1.01 & 4.01 & 0.29 & 8.79 \\
1.8 & 0.21 & 29.40 & 0.36 & 28.80 & 0.10 & 7.00 & 0.30 & 7.28 & 1.10 & 3.97 & 0.29 & 8.66 \\
1.9 & 0.37 & 27.72 & 0.46 & 26.91 & 0.11 & 6.90 & 0.34 & 7.06 & 1.25 & 3.88 & 0.28 & 8.50 \\
\bottomrule
\end{tabular}
\begin{tablenotes}[flushleft]
\small
\item[] \textit{Note:} The thresholds $u_{25}$, $u_{50}$, and $u_{75}$ are the 25th, 50th, and 75th quantiles of the simulated volatility path in each replication. See the notes to Table~\ref{tab:mc_continuous} for the definitions of the remaining entries.
\end{tablenotes}
\end{threeparttable}
\end{table}

\begin{table}[htbp]
\centering
\small
\begin{threeparttable}
\caption{Relative efficiency of the stability-based rule}
\label{tab:efficiency_all}
\setlength{\tabcolsep}{1.6pt}
\renewcommand{\arraystretch}{1.0}
\begin{tabular}{c c c c c c c c}
\toprule
$\beta$ & $g(x)=x^2$ & $g(x)=x$ & $g(x)=\sqrt{x}$ & $g(x)=\log(x)$ & $g(x)=\mathbf1_{\{x\le u_{25}\}}$ & $g(x)=\mathbf1_{\{x\le u_{50}\}}$ & $g(x)=\mathbf1_{\{x\le u_{75}\}}$ \\
\midrule
1.1 & 1.02 & 1.04 & 0.95 & 0.85 & 0.90 & 0.85 & 0.34 \\
1.2 & 1.03 & 1.04 & 0.96 & 0.88 & 0.92 & 0.86 & 0.36 \\
1.3 & 1.03 & 1.04 & 0.97 & 0.90 & 0.94 & 0.88 & 0.39 \\
1.4 & 1.04 & 1.04 & 0.98 & 0.92 & 0.95 & 0.89 & 0.41 \\
1.5 & 1.05 & 1.03 & 1.00 & 0.94 & 0.97 & 0.91 & 0.42 \\
1.6 & 1.07 & 1.02 & 1.01 & 0.96 & 0.99 & 0.92 & 0.44 \\
1.7 & 1.09 & 1.02 & 1.02 & 0.97 & 1.00 & 0.94 & 0.46 \\
1.8 & 1.11 & 1.01 & 1.03 & 0.98 & 1.02 & 0.96 & 0.46 \\
1.9 & 1.17 & 1.02 & 1.03 & 0.96 & 1.03 & 0.98 & 0.46 \\
\bottomrule
\end{tabular}
\begin{tablenotes}[flushleft]
\small
\item[] \textit{Note:} Relative efficiency is defined as $\mathrm{MAD}(p^*)/\mathrm{MAD}_{\mathrm{algo}}$, where $p^*$ is the ex ante optimum of Tables~\ref{tab:mc_continuous} and ~\ref{tab:mc_indicator}. Values above one arise because the rule adapts $p$ per replication, whereas $p^*$ is constrained to a single value across all replications.
\end{tablenotes}
\end{threeparttable}
\end{table}

\section{Empirical application}
\label{sec:empirical}

In this section, we apply the proposed estimators to two crypto
assets, Bitcoin (BTC) and TRUMP, which differ markedly in liquidity,
trading intensity, and jump behavior. BTC is the largest and most actively traded cryptocurrency, while TRUMP is a newly issued token associated with U.S. President Donald Trump, introduced in mid-January 2025, shortly before our sample begins.

The data consist of 1-minute prices from Binance covering January~20 to 
December~31, 2025, with timestamps in Coordinated Universal Time 
(UTC). The common sample window ensures comparability across assets 
and covers episodes of both tranquil trading and heightened market 
stress. Prices are sampled on the regular 1-minute grid, and returns are
computed as 1-minute log-price differences.

As a preliminary step, we estimate the jump-activity index $\beta$
of each asset over the full sample using the bipower-variation-based
procedure of \citet{kolokolov2022estimating}, obtaining
$\widehat\beta=1.90$ for BTC, in line with the pure-jump evidence for
Bitcoin reported there, and $\widehat\beta=1.56$ for TRUMP, indicating
more pronounced jump dominance.\footnote{Several consistent estimators
of the jump-activity index have been proposed in the literature (e.g.,
\citealp{ait2009estimating,todorov2011limit,todorov2015jump}).} These
estimates serve as inputs for the analysis that follows.

\subsection{Occupational characteristics of volatility}
\label{subsec:occupation}

We first study how the volatility processes of the two
cryptocurrencies differ over the sample period, in average level, in
variability, and in the frequency of extreme episodes. The occupational measures introduced in
Section~\ref{sec:functionals} are designed for such questions, as they
summarize the entire volatility path through the time it spends at
different levels.

The spot volatility process $\sigma_t$ itself is not comparable across assets
with different activity indices, because the same value of $\sigma_t$
implies different return distributions at different values of $\beta$.
We therefore rescale the estimated spot volatility into a common unit,
the expected magnitude of a one-minute return. Specifically, for each asset, define
the \emph{normalized spot volatility}
\[
\widehat\Sigma_{t} \equiv c_{\widehat\beta}(1)\,\Delta_{\mathrm{emp}}^{1/\widehat\beta}\,\widehat\sigma_{n,t},
\]
where $\widehat\sigma_{n,t}$ is the piecewise-constant spot volatility
process from \eqref{eq:piecewise-process}, $\widehat\beta$ is the
estimated activity index,
$c_{\widehat\beta}(1)=2^{1-1/\widehat\beta}\,\Gamma(1-1/\widehat\beta)/\pi$
follows from \eqref{eq:normal-const}, and $\Delta_{\mathrm{emp}}$ is a
fixed reporting horizon, set to one minute. Throughout, $\widehat\Sigma_t$ is reported in basis
points per one-minute interval. 

We summarize the occupational distribution of the normalized spot
volatility through five characteristics, namely the occupational mean,
median, standard deviation, interquartile range, and upper tail. These are the occupational analogues of the classical summary statistics of a data sample, with the distribution of a random variable replaced by the occupational distribution. As
these empirical quantities appear only in this section, we suppress
their dependence on $n$ and $T$ in the notation. The occupational mean
and standard deviation are
\begin{align}
\label{eq:occ-moments}
\widehat{\mu}
\equiv
\frac{1}{T}\int_0^T \widehat\Sigma_{s}\,ds,
\qquad
\widehat{\mathrm{SD}}
\equiv
\bigg(\frac{1}{T}\int_0^T
\big(\widehat\Sigma_{s}-\widehat{\mu}\big)^2\,ds\bigg)^{1/2}.
\end{align}
Following \eqref{eq:empirical-occupation-time} and
\eqref{volatility-quantile-estimator}, the normalized occupation time
of $\widehat\Sigma_{t}$ and the associated quantiles are
\[
\widehat F(x) \equiv \frac{1}{T}\int_0^T \mathbf 1_{\{\widehat\Sigma_{s}\le x\}}\,ds,
\qquad
\widehat Q(\alpha) \equiv \inf\{x:\widehat F(x)\ge\alpha \},
\]
which deliver the occupational median $\widehat Q(0.50)$, the
interquartile range $\widehat Q(0.75)-\widehat Q(0.25)$, and the upper
tail $\widehat Q(0.95)$.

{For each asset, the power index $p$ is obtained by applying the
stability-based rule of Section~\ref{subsec:algo} with indicator
targets $g_u(x)=\mathbf 1_{\{x\le u\}}$, where $u$ ranges over the
25th, 50th, and 75th sample quantiles of that asset's block-based spot
volatility estimates, and averaging the three resulting values. This
yields $p=0.61$ for BTC and $p=0.68$ for TRUMP.}\footnote{{Results
are virtually unchanged for nearby choices of each asset's power index
in $[0.50,0.75]$. The cross-asset ranking and the qualitative
conclusions in this section are unaffected.}}
{The same asset-specific values of $p$ are used throughout the
empirical application.}

\begin{figure}[h!]
\centering
\includegraphics[width=0.8\textwidth]{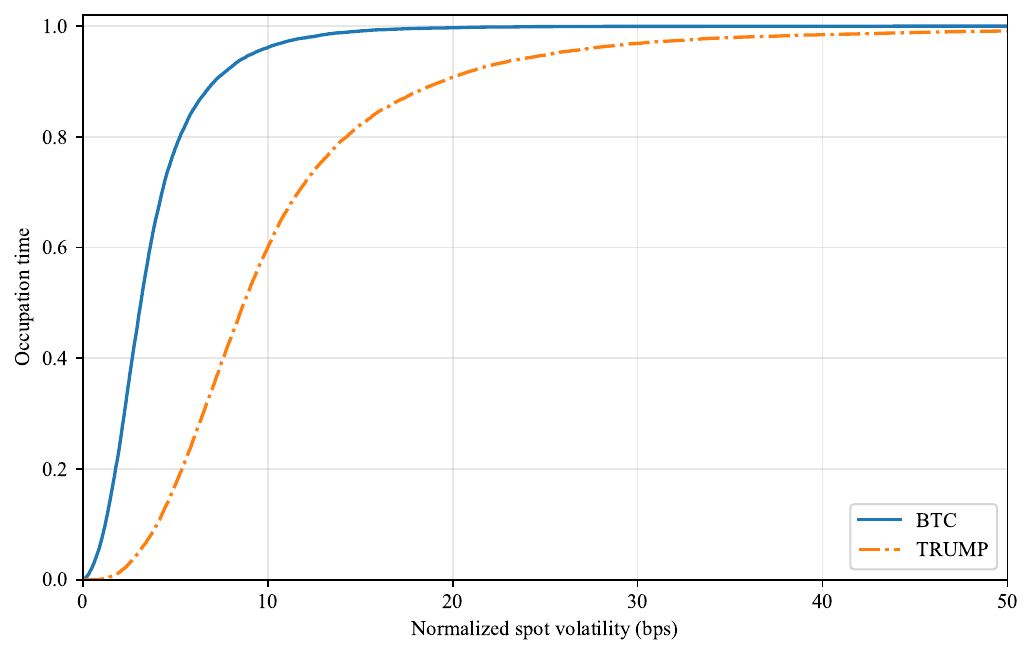}
\caption{{Occupation times $\widehat F(x)$ of the normalized spot
volatility for BTC and TRUMP, January~20 to
December~31, 2025, computed with asset-specific power indices
$p=0.61$ (BTC) and $p=0.68$ (TRUMP).}}
\label{fig:occupation}
\end{figure}

Figure~\ref{fig:occupation} plots the normalized occupation time
$\widehat F(x)$ for each asset, which represents the fraction of the
sample period during which $\widehat\Sigma_{t}$ remains below $x$. The
two occupation times display a clear contrast in risk profiles. For 
BTC, $\widehat F(x)$ rises rapidly and approaches one at relatively 
low values of $x$, indicating that its spot volatility is concentrated 
at small magnitudes. For TRUMP, $\widehat F(x)$ lies far to the right 
and rises much more gradually, suggesting prolonged exposure to 
high-risk regimes and a heavier upper tail.

\begin{table}[h!]
\centering
\small
\begin{threeparttable}
\caption{Occupational characteristics of the normalized spot volatility}
\label{tab:occupation_summary}
\setlength{\tabcolsep}{12pt}
\renewcommand{\arraystretch}{1}
\begin{tabular}{lccccc}
\toprule
Asset 
& Mean (bps)
& Median (bps)
& SD (bps)
& IQR (bps)
& UT (bps) \\
\midrule
BTC   & 3.81 & 3.13 & 2.94 & 2.68 & 9.18 \\
TRUMP & 11.02 & 8.71 & 12.29 & 6.85 & 25.36 \\
\bottomrule
\end{tabular}
\begin{tablenotes}[flushleft]
\small
\item[] \textit{Note:} Occupational summaries of $\widehat\Sigma_{t}$ in basis points. Mean and SD are defined in \eqref{eq:occ-moments}. Median $=\widehat Q(0.50)$, IQR $=\widehat Q(0.75)-\widehat Q(0.25)$, and UT (upper tail) $=\widehat Q(0.95)$ are computed from the occupation time shown in Figure~\ref{fig:occupation}. Each asset uses its own power index ($p=0.61$ for BTC and $p=0.68$ for TRUMP).
\end{tablenotes}
\end{threeparttable}
\end{table}

Table~\ref{tab:occupation_summary} quantifies this contrast through
the five occupational characteristics, which we read in three groups. The mean and the median
measure the average level of volatility. Both are roughly three times
larger for TRUMP, at 11.02~bps and 8.71~bps, than for BTC, at 3.81~bps
and 3.13~bps, and for each asset the mean exceeds the median, 
consistent with a right-skewed occupational distribution. The standard 
deviation and the interquartile range measure
the variation of volatility. Here the contrast sharpens, as TRUMP's SD of 12.29~bps is
about four times BTC's 2.94~bps, so TRUMP's spot volatility is not only
higher on average but also substantially more variable. Finally, the 
upper-tail quantile summarizes extreme volatility episodes. Spot
volatility stays below 9.18~bps for BTC during 95\% of the sample
period, whereas the corresponding level for TRUMP is 25.36~bps, nearly
three times higher.

\subsection{Occupational correlation and the leverage effect}
\label{subsec:leverage}

Beyond the univariate occupational summary statistics in
Section~\ref{subsec:occupation}, we
now study the bivariate dependence between the log-price level and the
normalized spot volatility. In equity markets, price declines tend
to be accompanied by rising volatility. \citet{black1976studies} and
\citet{christie1982stochastic} attributed this pattern to financial
leverage, as a fall in the stock price raises the firm's
debt-to-equity ratio and thereby the riskiness of its equity, and the
phenomenon has since been known as the leverage effect. Whether similar
dependence should be expected in crypto markets is unclear a priori,
because a token is not a claim on the cash flows of a levered firm and
the balance-sheet mechanism does not apply.

Measurement is a further issue in our setting. High-frequency work
on the leverage effect (e.g.,
\citealp{bandi2012time,wang2014estimation,ait2017estimation,kalnina2017nonparametric})
mostly targets the instantaneous correlation between the Brownian
shocks driving the price and its volatility, a quantity that is
degenerate under a pure-jump specification in which price increments
carry no Brownian component. The occupational correlation introduced in
Example~\ref{exp:correlation} offers an alternative measure of
price--volatility dependence that remains meaningful in this setting,
as it compares the levels of the two paths over the sample period
rather than their instantaneous shocks. This level-based view is also
closer to the evidence in the original studies, which documented the
co-movement between observed prices and estimated volatility over
calendar periods rather than an instantaneous correlation between
shocks, a quantity that would have been hard to discern by market observers. Since the
existing high-frequency literature reserves the term leverage effect for the
instantaneous notion, we adopt the separate name occupational
correlation for our measure and use it to examine the price--volatility
dependence of the two assets.

Following that example, we take the observable state variable to be the
log-price itself, $Y_t\equiv X_t$. Let
$\bar X\equiv T^{-1}\int_0^T X_{s}\,ds$ and
$\mathrm{SD}_X\equiv(T^{-1}\int_0^T(X_{s}-\bar X)^2\,ds)^{1/2}$ denote
the time-averaged mean and standard deviation of $X_t$. The occupational
covariance and correlation between $\widehat\Sigma_{t}$ and
$X_t$ are
\[
\widehat\Cov\equiv \frac{1}{T}\int_0^T
\bigl(\widehat\Sigma_{s}-\widehat\mu\bigr)
\bigl(X_{s}-\bar X\bigr)\,ds,
\qquad
\widehat\Corr\equiv
\frac{\widehat\Cov}{\widehat{\mathrm{SD}}\cdot\mathrm{SD}_X},
\]
where $\widehat\mu$ and $\widehat{\mathrm{SD}}$
are from \eqref{eq:occ-moments}. Since correlation is invariant
to scale, the normalization of the spot volatility is immaterial here,
and replacing $\widehat\Sigma_{t}$ with $\widehat\sigma_{n,t}$ yields
the same value of $\widehat\Corr$.

\begin{figure}[h!]
\centering
\includegraphics[width=1\textwidth]{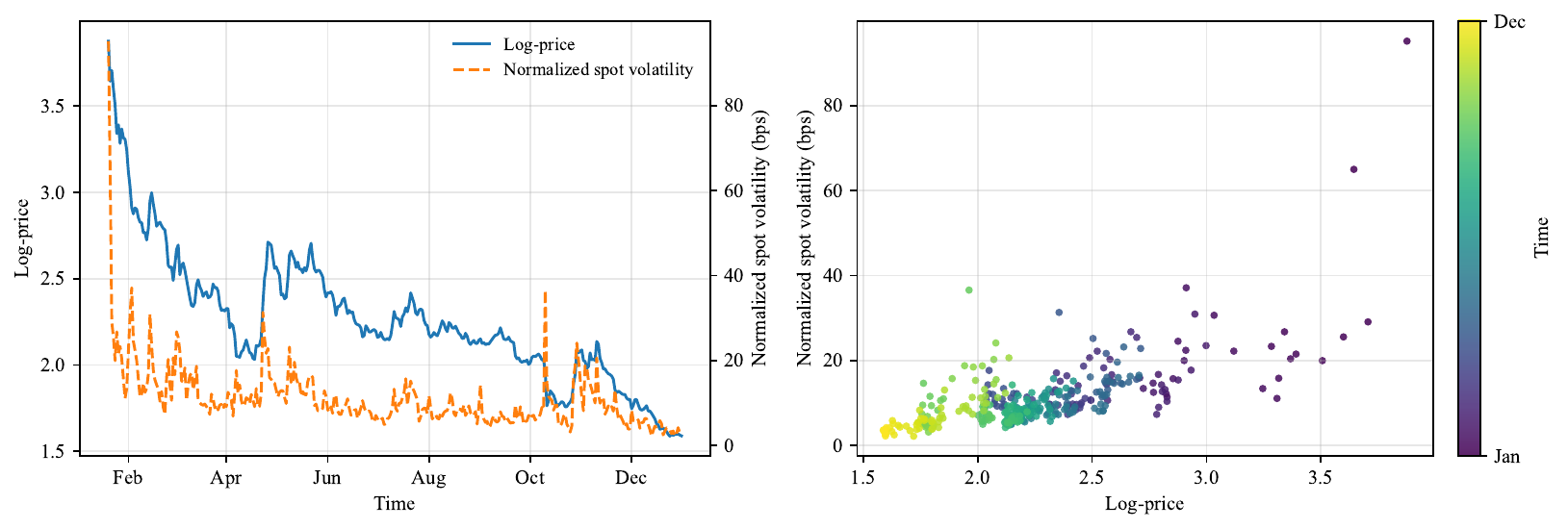}
\caption{Mechanism behind TRUMP's positive occupational
correlation. Left panel: daily averages of the log-price (left axis)
and the normalized spot volatility (right axis, basis points). Right
panel: daily means of log-price and normalized spot volatility, with
marker color recording calendar time from January to December.}
\label{fig:trump_leverage}
\end{figure}

The estimated occupational correlation is $-29.00\%$ for BTC and
$41.47\%$ for TRUMP, computed at the same asset-specific power indices
as in Section~\ref{subsec:occupation}.
The negative
estimate for BTC is consistent with leverage-type dependence in which
lower price levels coincide with higher spot volatility. The
positive estimate for TRUMP reflects the joint evolution of its price
and volatility after the token's launch. Over the sample, the price
declines from an early peak near \$48 to about \$4.80 by the end of
the year, and spot volatility subsides from its elevated initial
level over the same period. Deviations of the two paths from their
time averages therefore tend to share the same sign, which the
occupational correlation records as positive dependence.
Figure~\ref{fig:trump_leverage} illustrates this pattern.
The left panel shows the common downward trend in the daily
averages of the log-price and the normalized spot volatility, and the
right panel shows that days early in the sample combine high price
levels with high volatility while later days combine low levels of
both. The contrast between the
two assets indicates that the negative price--volatility dependence
documented for equities is not universal in crypto markets, consistent
with the absence of the balance-sheet mechanism underlying the
original leverage hypothesis.

\subsection{Volatility around information arrivals}
\label{subsec:events} 

We conclude the empirical analysis by examining how spot volatility
behaves around discrete information arrivals. We consider two observable
sources of arrivals, namely monetary policy announcements and posts by
President Trump on Truth Social that are classified as financially
relevant. The two sources differ in nature, as the former carries
market-wide news while the latter is tied to a single originator, and
we therefore analyze them separately for both assets.

Monetary policy announcements are taken from the official FOMC
calendar and cover scheduled meetings, statements, minutes releases,
press conferences, speeches and testimonies, Beige Book publications,
and major macroeconomic data releases. Each item carries an exact
release time, which we record in Coordinated Universal Time (UTC) and
align with the high-frequency price data.

Each Truth Social post by President Trump during the sample period is 
treated as a candidate information arrival. To 
distinguish economically relevant information from non-financial 
communication, we classify each post with an LLM, labeling it 
financially relevant if its content concerns cryptocurrencies, 
financial markets, or macroeconomic policy. The 
LLM is used solely as a deterministic text-based classification tool, 
with no access to price or volatility data. The prompt, decoding 
parameters, and random seeds are fixed ex ante, ensuring that the 
classification cannot adapt to realized outcomes. Robustness to model 
choice is verified by replicating the classification across multiple 
open-source LLMs (\texttt{llama3.2:3b}, \texttt{deepseek-r1:8b}, 
\texttt{deepseek-r1:14b}, \texttt{qwen3:8b}). The resulting labels and
all empirical conclusions are essentially unchanged. See 
Appendix~B in the supplementary materials for details.

For each arrival, we form a symmetric event window extending one
hour on each side of the arrival time. This contrasts near-arrival
periods with ordinary times, rather than a local post- versus
pre-announcement comparison for which the pre-announcement hour is a
poor control when news is anticipated. Within each source, overlapping
windows are merged, event time is the union of the windows, and the
remaining periods constitute non-event time. We then compare the
occupation times of the normalized spot volatility over event and
non-event periods, computed at the same asset-specific power
indices as in Section~\ref{subsec:occupation}.\footnote{Results are qualitatively
robust to alternative window lengths of 30, 90, and 120 minutes.
Details are omitted for brevity.}

\begin{figure}[h!]
\centering
\includegraphics[width=1\textwidth]{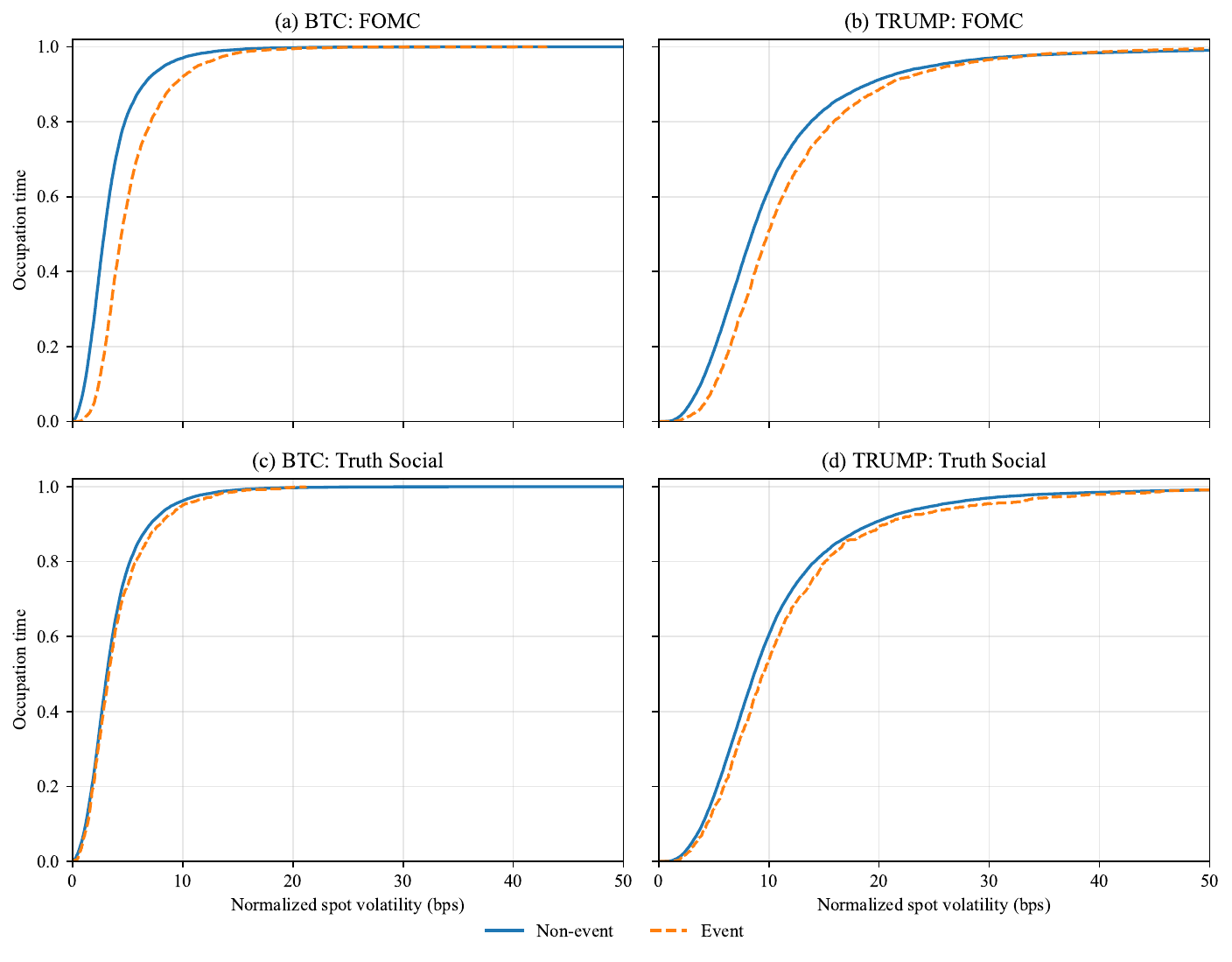}
\caption{Occupation times of the normalized spot volatility for BTC
and TRUMP over event and non-event periods. Panels~(a) and~(b) use
FOMC-related announcements as events for BTC and TRUMP, respectively,
and panels~(c) and~(d) use financially relevant Truth Social posts.
Event windows are $\pm 60$-minute intervals around the arrivals, with
overlapping windows merged, and non-event periods are the complement.}
\label{fig:event_occupation}
\end{figure}

Figure~\ref{fig:event_occupation} reports the results, with the
two assets in columns and the two event sources in rows. In every
panel, the occupation time for event windows lies to the right of that
for non-event periods, consistent with first-order stochastic
dominance, so spot volatility is elevated around information
arrivals. The magnitude of the shift differs sharply
across sources. FOMC-related announcements shift the occupational
distribution visibly for both BTC and TRUMP, consistent with
market-wide policy news that raises spot volatility across crypto
assets. Truth Social posts produce a much more muted shift,
discernible mainly for TRUMP. This pattern is natural if the posts mix
high-signal and low-signal content, as aggregation across
heterogeneous messages attenuates the average effect, while the
TRUMP-specific displacement points to a token-linked attention channel
that is weaker for BTC.\footnote{For example, on March~23, 2025, President
Trump posted ``I LOVE \$TRUMP --- SO COOL!!! The Greatest of them
all!!!!!!!!!!!!!!!!''
(\url{https://truthsocial.com/@realDonaldTrump/posts/114212473050636851}), which is specific to the TRUMP token.} We therefore read the
Truth Social evidence as reflecting asset-specific attention rather
than a market-wide volatility shock of the FOMC type.\footnote{The
comparison is descriptive rather than causal, since elevated
volatility around events may reflect news content, attention
spillovers, or correlated trader behavior, and our analysis is not
designed to disentangle these channels.}

\section{Concluding remarks}
\label{sec:conclusion}

This paper has developed a uniform estimation theory for
integrated volatility functionals when the underlying price is a
pure-jump semimartingale, covering volatility occupation times and
quantiles as well as plug-in $M$-estimation, and it has proposed a
stability-based selection rule that makes the procedure operational. We demonstrate the empirical relevance of the proposed estimators in a study of cryptocurrencies.

The natural next step is the associated distribution theory. A
central limit theory for the empirical occupation measure would turn
the consistency results obtained here into inference procedures for
occupation times, quantiles, and $M$-estimators. The stable regime
poses distinctive challenges for this program, as the Gaussian limit
arguments available in the Brownian setting no longer apply. This important extension is beyond the scope of the current paper and is left for future research.

\spacingset{1}
\bibliographystyle{apalike}
\bibliography{reference}

\end{document}